\documentclass[a4paper]{amsart}
\usepackage{epsf}
\usepackage[all]{xy}
\newcommand{\conf}{{\em cf.\ }}

\newtheorem*{theorem}{Theorem}
\newtheorem*{lemma}{Lemma}
\newtheorem*{proposition}{Proposition}
\newtheorem*{corollary}{Corollary}

\theoremstyle{definition}

\theoremstyle{remark}

\numberwithin{equation}{section}

\newcommand{\mysection}[2]{\subsection{#1}\label{#2}}

\newcommand{\internalcomment}[1]{}

\newcommand{\Z}{\mathbf{Z}}

\newcommand{\ko}{\: , \;}

\renewcommand{\tilde}[1]{\widetilde{#1}}

\newcommand{\ra}{\rightarrow}
\newcommand{\la}{\leftarrow}
\newcommand{\da}{\downarrow}

\newcommand{\arr}[1]{\stackrel{#1}{\rightarrow}}
\newcommand{\longarr}[1]{\xrightarrow{#1}}

\newcommand{\iso}{\stackrel{_\sim}{\rightarrow}}
\newcommand{\liso}{\stackrel{_\sim}{\leftarrow}}

\newcommand{\opname}[1]{\operatorname{\mathsf{#1}}}

\newcommand{\Mod}{\opname{Mod}\nolimits}
\newcommand{\Comod}{\opname{Comod}\nolimits}

\newcommand{\cdg}{\opname{cdg}\nolimits}
\newcommand{\id}{\mathbf{1}}
\newcommand{\R}{\mathbf{R}}

\newcommand{\ten}{\otimes}
\newcommand{\lten}{\otimes^{\mathsf{L}}}
\newcommand{\lrelten}{\otimes^{\mathsf{L},\mbox{\scriptsize rel}}}
\newcommand{\tp}[1]{^{\ten #1}}

\renewcommand{\ker}{\opname{ker}\nolimits}

\newcommand{\op}{^{op}}
\newcommand{\Zy}[1]{\opname{Z}^{#1}}

\renewcommand{\H}[1]{{H}^{#1}}

\newcommand{\Hs}{{H}^*}
\newcommand{\HH}[1]{{HH}^{#1}\,}

\newcommand{\HHs}{{HH}^*}
\newcommand{\HHsp}{{HH}^{*+1}}

\newcommand{\cd}{{\mathcal D}}

\newcommand{\cf}{{\mathcal F}}

\newcommand{\cw}{{\mathcal W}}

\newcommand{\eps}{\varepsilon}

\renewcommand{\phi}{\varphi}

\newcommand{\Hom}{\opname{Hom}}

\newcommand{\Ext}{\opname{Ext}}

\newcommand{\Aut}{\opname{Aut}}

\newcommand{\Der}{\opname{Der}}
\newcommand{\Coder}{\opname{Coder}}
\newcommand{\Lie}{\opname{Lie}}
\newcommand{\dpic}{\opname{DPic}}
\newcommand{\liedpic}{\opname{LieDPic}}
\newcommand{\lieaut}{\opname{LieAut}}
\newcommand{\defo}{\opname{Defo}}

\newcommand{\n}{\mathfrak{n}}
\newcommand{\m}{\mathfrak{m}}

\begin{document}

\title{Hochschild cohomology and derived Picard groups}

\author{Bernhard Keller}
\address{UFR de Math\'ematiques\\
   UMR 7586 du CNRS \\
   Case 7012\\
   Universit\'e Paris 7\\
   2, place Jussieu\\
   75251 Paris Cedex 05\\
   France }

%\thanks{Partially supported by the European Network `Algebraic
%Lie Representations', Contract ERB-FMRX-CT97-0100}

\email{
\begin{minipage}[t]{5cm}
keller@math.jussieu.fr \\
www.math.jussieu.fr/ $\tilde{ }$ keller
\end{minipage}
}

\date{January 2002, last modified on September 25, 2003}

\dedicatory{Dedicated to Idun Reiten on the occasion of her
sixtieth birthday}

\begin{abstract}
We interpret Hochschild cohomology as the Lie algebra
of the derived Picard group and deduce that it is preserved
under derived equivalences.
\end{abstract}

\maketitle

%\tableofcontents

\section{Introduction}

The Hochschild cohomology groups $\HH{i}(A,A)$ of an algebra $A$
over a field $k$ can be interpreted as higher extension groups of
the bimodule $A$ by itself or as morphisms from $A$ to $A[i]$ in
the derived category $\cd(A\op\ten A)$ of $A$-$A$-bimodules. This
last interpretation shows that they are preserved under derived
equivalences \cite{Rickard89}, i.e. if $X$ is a complex of
$A$-$B$-bimodules such that the total derived tensor product by
$X$ is an equivalence $\cd A \ra \cd B$, then $X$ yields a natural
isomorphism from $\HH{i}(A,A)$ to $\HH{i}(B,B)$. This isomorphism
is compatible with the cup product, since the cup product
corresponds to the composition of morphisms in the derived
category of bimodules. However, it is not clear whether the
isomorphism given by $X$ also respects the Gerstenhaber bracket on
Hochschild cohomology \cite{Gerstenhaber63}. We will show that
this is indeed the case by providing an intrinsic interpretation
of the Gerstenhaber bracket in terms of derived categories. The
basic idea is to view Hochschild cohomology as an analogue of the
Lie algebra associated with an algebraic group (more precisely, a
group-valued functor). This group will be the derived Picard group
\cite{Zimmermann96} \cite{RouquierZimmermann98} \cite{Yekutieli98}
of $A$ (more precisely, the functor which sends a commutative
differential graded $k$-algebra $R$ to the $R$-relative derived
Picard group of $A\ten_k R$). Our interpretation generalizes the
fact that the first Hochschild cohomology group of a
finite-dimensional algebra $A$ is the Lie algebra of the group of
outer automorphisms of $A$.

The author thanks the referee for his careful reading of the manuscript and
for many helpful remarks.

\section{Reminder on derived equivalences}

\mysection{Derived categories and the Hochschild cohomology algebra}{DerCatHochAlg}
Let $k$ be a field and $A$ a $k$-algebra, i.e. an associative unital $k$-algebra.
Let $\Mod A$ denote the category of right $A$-modules. Let $\cd A$ denote the
(unbounded) derived category of $\Mod A$. Thus, the objects of $\cd A$ are all
complexes
\[
\ldots \ra M^p \arr{d} M^{p+1} \ra \ldots \ko p\in \Z\ko d^2=0\ko
\]
of right $A$-modules and its morphisms are obtained from morphisms of complexes
by formally inverting all quasi-isomorphisms, i.e. morphisms of complexes inducing
isomorphisms in homology. Let us recall the most basic examples of morphisms
in the derived category: We identify an $A$-module $L$ with the complex whose
$0$-component is $L$ and whose components in all other degrees vanish. Then,
if $L$ and $M$ are $A$-modules, the group of morphisms in $\cd A$ from $L$
to $M$ identifies with the group of $A$-linear maps from $L$ to $M$ and,
more generally, we have a natural isomorphism
\[
\Hom_{\cd A}(L,M[i]) \iso \Ext_A^i(L,M)
\]
for each $i\in\Z$, where, for a complex $K$, we denote by $K[i]$ the complex
with components $K[i]^p=K^{i+p}$ and differential $(-1)^i d_K$. By convention,
$\Ext$-groups vanish in negative degrees.

In particular, we can identify the Hochschild cohomology groups $\HHs(A,A)$ with groups
of morphisms in the derived category of $A$-$A$-bimodules: Indeed, since $A$
is flat over $k$, we have a canonical isomorphism
\[
\HH{i}(A,A) \iso \Ext^i_{A\op\ten A}(A,A)
\]
and thus a canonical isomorphism
\[
\HH{i}(A,A) \iso \Hom_{\cd (A\op\ten A)}(A,A[i])\ko i\in \Z.
\]
Under this isomorphism, the cup product on Hochschild cohomology corresponds to
the graded composition in the derived category. More precisely, the cup product
of the cohomology classes corresponding to $f: A \ra A[j]$ and $g: A \ra A[i]$
corresponds to the composition $f[i] \,\circ\, g$.

\mysection{Derived equivalence}{DerEq}
Let $A$ and $B$ be two $k$-algebras. We recall one version of Rickard's Morita
theorem for derived categories \cite{Rickard89} \cite{Rickard91}.

\begin{theorem} The following are equivalent
\begin{itemize}
\item[(i)] There is a triangle equivalence $F: \cd A \iso \cd B$.
\item[(ii)] There are bimodule complexes $X\in \cd(A\op\ten B)$ and $Y\in \cd(B\op\ten A)$
and isomorphisms
\[
X\lten_B Y \iso A \mbox{ in } \cd(A\op\ten A) \mbox{ and } Y \lten_A X \iso B \mbox{ in } \cd(B\op\ten B).
\]
\end{itemize}
\end{theorem}

Here the symbol $\lten$ denotes the total derived tensor functor
\cite{Spaltenstein88}.
The implication from (ii) to (i) is easy: Indeed, the functor
\[
F=\, ?\lten_A X\, : \cd A \ra \cd B
\]
is an equivalence whose inverse is given by $?\lten_B Y$, \conf \cite{Rickard91}.
The implication from (ii) to (i) is considerably more delicate. One can also show
\cite{Rickard91} that if
$X\in \cd(A\op\ten B)$ is a bimodule complex such that the associated
functor $?\lten_A X : \cd A \ra \cd B$ is an equivalence, then (ii) holds
for $X$ and
\[
Y= \mathsf{R}\Hom_B(X,B).
\]
Thus the essential datum is that of $X$. We call such $X$ an {\em invertible}
bimodule complex and $Y$ its {\em inverse}.
Two algebras $A$ and $B$ are
called {\em derived equivalent} if the conditions of the theorem hold.
If we consider other variants of the derived categories (e.g. the
bounded derived categories), we obtain the same equivalence relation
on the class of $k$-algebras, \conf \cite{Rickard89}. Of course, derived equivalence
generalizes Morita equivalence. The following example, a particular case of Koszul duality
\cite{BeilinsonGinzburgSoergel96}, \cite{GreenReitenSolberg2000},
shows that this generalization is non trivial.

\mysection{An example}{Example}
Let $V$ be a vector space of dimension $n+1$, denote by $S^i$ the $i$-th
symmetric power of its dual space and by $\Lambda^i$ the $i$-th exterior
power of $V$. Let $A$ be the algebra of upper triangular matrices
\[
A= \left( \begin{array}{cccc}
S^0 & S^1 & \ldots & S^n \\
0   & S^0 & \ldots & S^{n-1} \\
\vdots & \vdots & \ldots  & \vdots \\
0   & 0   & \ldots & S^0
\end{array} \right)
\]
and $B$ the algebra of lower triangular matrices
\[
B=\left( \begin{array}{cccc}
\Lambda^0 & 0 & \ldots & 0 \\
\Lambda^1 & \Lambda^0 & \ldots & 0 \\
\vdots    & \vdots    & \ldots & \vdots \\
\Lambda^n & \Lambda^{n-1} & \ldots & \Lambda^0
\end{array} \right) .
\]
Then $A$ is derived equivalent to $B$, but for $n>2$, $A$ is not
Morita equivalent to $B$. If $E_i$ denotes the simple $A$-module
corresponding to the projection on the $(i+1)$-th diagonal component,
$i=0,\ldots, n$, then the complex
\[
\ldots \ra 0 \ra E_n \arr{0} E_{n-1} \arr{0} \ldots \arr{0} E_{0} \ra 0 \ldots
\]
is the restriction to $A$ of a bimodule
complex $X\in \cd(A\op\ten B)$, unique up  to isomorphism \cite{Keller93},
whose associated tensor functor is an equivalence $\cd A \ra \cd B$.
Historically, this example comes from geometry: Beilinson showed in
\cite{Beilinson78} that the derived categories
of both $A$ and $B$ are triangle equivalent to the derived category
of the category of quasicoherent sheaves on the projective space $\mathbf{P}(V)$.

\mysection{Invariance of the algebra structure}{InvAlg}
The following theorem is due to Dieter Happel \cite{Happel89} in the special
case of derived equivalences coming from tilting modules and to
Jeremy Rickard \cite{Rickard91} in the general case. Let $A$ and $B$ be
derived equivalent algebras and $X\in \cd(A\op\ten B)$ and invertible
bimodule complex

\begin{theorem} There is a canonical algebra isomorphism
\[
\phi_X : \HHs(A,A) \ra \HHs(B,B).
\]
\end{theorem}

After the preparations we have made, it is easy to construct $\phi_X$:
Indeed, let $Y$ be the inverse of $X$ with isomorphisms
$u: Y\lten_A X \iso B$ and $v: X \lten_B Y \iso A$.
Let $\phi_{X,u}$ send $f:A \ra A[i]$ to
\[
u[i] \circ (Y\lten_A f \lten_A X) \circ u^{-1}\; : \; B \ra B[i].
\]
If $u'$ is another isomorphism from $Y\lten_A X$ to $B$, then
$u'=zu$ for an invertible central element $z$ of $B$. So we have
\[
\phi_{X,u'}(f) = z[i]\, \phi_{X,u}(f) \, z^{-1}
\]
and this equals $\phi_{X,u}(f)$ since the center of $B$ is central
in $\HHs(B,B)$. We define $\phi_X=\phi_{X,u}$.

%%%%%%%%%%%%%%%%%%%%%%%%%%%%%%%%%%%%%%%%%%%%%%%%%%%%%%%%%%%%%%%%%%%%%%%%%%%%%%%%%%%%%%

\section{Reminder on the Gerstenhaber bracket}
Let $A$ be an algebra. The {\em Hochschild complex} of $A$ is the complex
\[
C^\cdot(A,A)= (A \ra \Hom_k(A,A) \ra \ldots \ra \Hom_k(A\tp{p}, A) \ra \ldots )
\]
with $A$ in degree $0$ whose differential maps $a\in A$ to $[a,?]$, and,
more generally, a $p$-cochain $c\in \Hom_k(A\tp{p},A)$ to the
$(p+1)$-cochain $dc$ defined by
\begin{eqnarray*}
(dc)(a_0,\ldots, a_p) & = & a_0 c(a_1,\ldots, a_p) - c(a_0 a_1, a_2,\ldots, a_p) + \ldots  \\
                      &   & + (-1)^{p+1} c(a_0, \ldots, a_{p-1}) a_p.
\end{eqnarray*}
The homology in degree $i$ of $C^\cdot(A,A)$ is $\HH{i}(A,A)$. For
a $p$-cochain $c_1$, a $q$-cochain $c_2$ and an integer $0\leq i \leq p-1$, define
a $(p+q-1)$-cochain $c_1\bullet_i c_2$ by
\[
(c_1 \bullet_i c_2)(a_1, \ldots, a_{p+q-1}) =
c_1(a_1, \ldots, a_i, c_2(a_{i+1}, \ldots, a_{i+q}), a_{i+q+1}, \ldots, a_{p+q-1}).
\]
Then the {\em Gerstenhaber product} is defined by
\[
(c_1 \bullet c_2)= \sum_{i=0}^{p-1} (-1)^{i(q-1)} c_1 \bullet_i c_2.
\]
%The terms on the right hand side are conveniently symbolized by the graphs
%\[
%\]
The Gerstenhaber product is not associative. However, the associator
\[
A(c_1, c_2, c_3)= (c_1 \bullet c_2)\bullet c_3 - c_1\bullet (c_2 \bullet c_3)
\]
is super symmetric in $c_2$ and $c_3$ endowed with suitable degrees (super=$\Z/2\Z$-graded).
Namely, we have
%Indeed, if we compute a composition
%$(c_1 \bullet c_2)\bullet c_3$, we find terms corresponding to the three following types
%of graphs
%\[
%\]
%The terms in the middle are exactly those occurring in $c_1\bullet (c_2 \bullet c_3)$.
%We see that the set of terms of the associator
%\[
%A(c_1, c_2, c_3)= (c_1 \bullet c_2)\bullet c_3 - c_1\bullet (c_2 \bullet c_3)
%\]
%is symmetric in $c_2$ and $c_3$. A more careful analysis taking the signs into account
%reveals that
\[
A(c_1, c_2, c_3) = (-1)^{(q-1)(r-1)} \, A(c_1, c_3, c_2)\ko
\]
where $c_2$ is a $q$-cochain and $c_3$ an $r$-cochain. Therefore the super commutator of the
Gerstenhaber bracket behaves like the commutator of an associative product, i.e. it
satisfies the super Jacobi identity. More precisely, we have the

\begin{lemma}\cite{Gerstenhaber63}
Endowed with the {\em Gerstenhaber bracket} defined by
\[
[c_1, c_2] = c_1 \bullet c_2 - (-1)^{(p-1)(q-1)} c_2 \bullet c_1\ko c_1 \in C^p(A,A)\ko c_2\in C^q(A,A)\ko
\]
and the differential $-d$, the graded space $C^{\cdot+1}(A,A)$ becomes a differential graded
Lie algebra. In particular, the homology $\HHsp(A,A)$ becomes a $\Z$-graded super Lie algebra.
\end{lemma}

For example, let $V$ be a vector space and $A$ the algebra of polynomial functions on $V$.
Then $\HH{1}(A,A)$ identifies with the space $Der_k(A,A)$ of $k$-linear derivations of $A$
and $\HH{i+1}(A,A)$ is canonically isomorphic to the exterior power $\Lambda^{i+1}_A \Der_k(A,A)$. Under
this isomorphism, the Gerstenhaber bracket corresponds to the Nijenhuis-Schouten-bracket,
which is the natural extension of the commutator of derivations.

\section{The Lie algebra of the derived Picard group}

\label{MainResult}
We will interpret Hochschild cohomology with the Gerstenhaber bracket as the Lie algebra
of a `generalized algebraic group', namely a group valued functor defined on a category
of commutative algebras. For this, let us recall the construction of the Lie algebra
of an algebraic group: Let $G$ be an algebraic group over $k$ considered as a
group-valued functor
\[
G : \{\mbox{commutative $k$-algebras}\} \ra \{ \mbox{groups} \} \ko R \mapsto G(R).
\]
Then the Lie algebra of $G$ is the space of tangent vectors at the origin, i.e.
\[
\Lie(G) = \ker( G(k[\eps]/(\eps^2)) \ra G(k) ).
\]
The bracket is induced by the commutator in $G$. To make this last statement more intuitive,
consider the example where $G=GL_n$. We have
\[
\Lie(GL_n) = \{ 1+ \eps X \; | \; X \in M_n(k) \; \}
\]
and the Lie bracket is determined by the identity
\[
(1+\eps_1 X_1)(1+ \eps_2 X_2)(1+\eps_1 X_1)^{-1} (1+\eps_2 X_2)^{-1} = 1+ \eps_1 \eps_2 \; [X_1, X_2]
\]
in $GL_n(k[\eps_1,\eps_2]/(\eps_1^2,\eps_2^2))$.

We will now define a group valued functor $\dpic_A$ whose Lie algebra will be
the Hochschild cohomology of $A$. Since this Lie algebra is graded, the category on
which $\dpic_A$ is defined should include the category of graded commutative algebras.
It turns out that a reasonable category is $\cdg k$, the category of commutative
differential graded $k$-algebras. To define
\[
\dpic_A : \cdg k \ra \{\mbox{groups}\}\ko
\]
we need the relative derived category (\conf \cite[Sect. 7]{Keller98}):
Let $R$ be a commutative differential graded
algebra (for example the algebra $k[\eps]/(\eps^2)$, where $\eps$ has any integer degree
and $d=0$). Let $E$ be a (typically non commutative) differential graded $R$-algebra.
The {\em $R$-relative derived category} $\cd_R E$ has as objects all differential
graded $E$-modules (these are precisely the complexes of $E$-modules if $E$ is
concentrated in degree $0$). The morphisms of $\cd_R E$ are obtained from morphisms
of differential graded $E$-modules by formally inverting all {\em $R$-relative quasi-isomorphisms},
i.e. all morphisms $s:L \ra M$ of differential graded $E$-modules whose restriction to
$R$ is an homotopy equivalence. For example, the relative derived category $\cd_k E$
equals the usual derived category $\cd E$ of the differential graded algebra $E$.

Rouquier-Zimmermann \cite{Zimmermann96} \cite{RouquierZimmermann98} and Yekutieli \cite{Yekutieli98} have
independently defined the derived Picard group of a ring. We generalize this as follows:
Let $R$ be a commutative differential graded algebra and $A$ an algebra. A bimodule
complex $U\in \cd_R(R\ten A\op\ten A)$ is {\em $R$-semifree}
if its underlying graded $R$-module is free; it is {\em invertible} if it is
$R$-semifree\footnote{presumably, all semifreeness conditions are redundant.}
and there exists an $R$-semifree bimodule complex $V\in \cd_R(R\ten A\op\ten A)$
such that there are isomorphisms
\[
U\lrelten_{A\ten R} V \iso R\ten A \mbox{ and } V\lrelten_{A\ten R} U \iso R\ten A
\]
in $\cd_R(R\ten A\op\ten A)$. The {\em $R$-relative derived Picard group of $A$} is
the set of isomorphism classes of invertible bimodule complexes $U$ in
$\cd_R(R\ten A\op\ten A)$. This set is endowed with the group law induced by the
derived tensor product. This group is denoted by $\dpic_A(R)$. It is functorial
with respect to $R$ so that we do obtain a functor
\[
\cdg k \ra \{\mbox{groups}\}.
\]
If $A$ is derived equivalent to $B$ and $X\in\cd(A\op\ten B)$ is
an invertible bimodule complex with inverse $Y$, then we have an isomorphism
\[
\dpic_A(R) \iso \dpic_B(R) \ko U \mapsto Y \lten_A U \lten_A X .
\]
In this sense, $\dpic_A$ is also functorial with respect to invertible bimodule
complexes $X\in \cd(A\op\ten B)$.

We now define the Lie algebra of $\dpic_A$.  Fix a degree $i\in \Z$
and let $R$ be the commutative differential graded algebra
$k[\eps]/(\eps^2)$, where $\eps$ is of degree $-i$ and $d=0$.  By
definition, $\liedpic_A^i$ is the set of isomorphism classes $U$ of
$\cd_R(R\ten A\op\ten A)$ such that $U$ is free as a graded $R$-module
and $U\ten_R k$ is isomorphic to $A$ in
$\cd(A\op\ten A)$.  The graded space $\liedpic_A^*$ is endowed with a
super Lie bracket defined as for algebraic groups
(\conf \ref{DefLieBracket}). For a super Lie
algebra $L$, we denote by $L\op$ the super Lie algebra with the
opposite bracket.  Its appearance in the theorem below is due to the
fact that we consider right modules.

\begin{theorem} There is a canonical isomorphism of graded super Lie algebras
\[
\HHsp(A,A)\op \iso \liedpic_A^*
\]
functorial with respect to invertible bimodule complexes $X\in \cd(A\op\ten B)$.
In particular, the Gerstenhaber bracket on $\HHsp(A,A)$ is preserved under
derived equivalence.
\end{theorem}

The rest of the article is devoted to the proof of the theorem. In fact,
we will prove it more generally for a differential graded $k$-algebra $A$.

\section{Proof of the main theorem}

\mysection{Outline of the proof}{OutlineProof} Let $A$ be a dg
algebra.  The idea is to construct an intermediate \lq
differential graded formal group $G$\rq\  whose Lie algebra is the
Hochschild complex with the Gerstenhaber bracket and which acts on
the relative derived category via bimodules. This group is the
group $G$ of automorphisms of the cobar construction $C^+$ of $A$,
where the cobar construction is viewed as a differential graded
counital (but not coaugmented) coalgebra (\conf
\ref{CoalgAutomBimod}). It naturally acts on the category of
differential graded comodules over $C^+$. Via the
bar-cobar-adjunction at the module level, this action translates
into an action of $G$ on the derived category via bimodules. On
the other hand, the Lie algebra of $G$ is the Lie algebra of
coderivations of $C^+$ and, by Stasheff's interpretation
\cite{Stasheff93}, this Lie algebra is the Hochschild complex
endowed with the Gerstenhaber bracket. This programme yields a Lie
algebra morphism
\[
\xymatrix{
\Zy{0}(C^{\cdot+1}(A,A)\ten \m) \ar[r] & \liedpic(A,R)
}
\]
for each augmented commutative dg algebra $R=k\oplus\m$ with $\m^2=0$.
To check that it induces an isomorphism, we need to identify
the set of deformation classes $\liedpic(A,R)$ with the group
\[
\Hom_{\cd(A\ten A\op)}(A,A\ten\m[1]).
\]
We will also need to know how the Lie algebra structure is reflected
under this identification. This is what we study first, in
the paragraphs \ref{InfDefModules} -- \ref{DefLieBracket} below.

\mysection{Infinitesimal deformations of modules}{InfDefModules}

Let $k$ be a field and $S$ a commutative dg $k$-algebra. Suppose that
$R$ is an augmented commutative dg $S$-algebra and denote by $\n$ the
kernel of the augmentation $R\ra S$. Thus we have the decomposition $R=S\oplus \n$.

Let $A$ be a (typically noncommutative) dg $S$-algebra, free as a graded $S$-module.
Then $A\ten_S R$ is a dg $R$-algebra. We consider the reduction functor
\[
\cd_R (A\ten_S R) \ra \cd_S A \ko L \mapsto L\ten_R S.
\]
We will study the fibers of this functor: Let $M$ be a dg $A$-module
which is free as a graded $S$-module.
Let $\cf$ be the category whose objects are the {\em deformations of $M$},
i.e. the pairs $(L,u)$ formed by a dg $A\ten R$-module $L$, free as a graded $R$-module,
and an isomorphism of $\cd_S A$
\[
u: L\ten_R S \iso M .
\]
Morphisms from $(L,u)$ to $(L', u')$ are given by morphisms
$v: L \ra L'$ of $\cd_R (A\ten_S R)$ such that $u' \circ (v\ten_R S) = u$.
We denote by
\[
\defo(M, R\ra S)
\]
the set of isomorphism classes of $\cf$. We denote by
\[
\defo'(M,R\ra S)
\]
the set of isomorphism classes of {\em weak deformations of $M$}, i.e.
dg $A\ten R$-modules $L$ free as graded $R$-modules such that $L\ten_R S$
is isomorphic to $M$. Note that we have an obvious forgetful map
\[
\defo(M,R\ra S) \ra \defo'(M,R\ra S).
\]
The group of automorphisms of $M$ in $\cd_S A$ acts on $\defo(M,R\ra S)$
via $(L,u).f=(L,f^{-1}\circ u)$ and the forgetful map clearly induces a bijection
\[
\defo(M,R\ra S)/ \Aut_{\cd_R A}(M) \iso \defo'(M,R\ra S).
\]

From now on and to the end of this section, we suppose that $\n^2=0$.
We will parametrize the connected components of $\cf$.
Let $(L,u)$ be an object of $\cf$. Since $L$ is free as a graded $R$-module,
the morphism of complexes $L\ten_R \n \ra L\n$ is invertible and the sequence
\[
0 \ra L\ten_R \n \ra L \ra L\ten_R S \ra 0
\]
is an exact sequence of dg $A$-modules which splits as a sequence of dg $S$-modules.
Thus it gives rise to a canonical triangle of $\cd_S A$ (but not of $\cd_R A$ !)
\begin{equation}
\label{CharTriangle}
L\ten_R \n \ra L \ra L\ten_R S \arr{\eps'} (L\ten_R \n)[1]
\end{equation}
Since $\n^2=0$, we have a canonical isomorphism
of dg modules
\[
L\ten_R \n \iso (L\ten_R S) \ten_S \n.
\]
Therefore, we can define a canonical morphism $\eps(L,u)$ of $\cd_S A$
by the commutative square
\[
\xymatrix{
L\ten_R S  \ar[r]^(0.45){\eps'} \ar[d]_{u} & L \ten_R \n[1] \ar[r]^(0.4){\sim} &
                                                (L\ten_R S) \ten_S \n[1] \ar[d]_{u\ten\n[1]} \\
M \ar[rr]^{\eps(L,u)} &  &  M\ten_S \n[1]
}
\]
Clearly the morphism $\eps(L,u)$ only depends on the isomorphism class of $(L,u)$
in the category $\cf$.

\begin{proposition}
The map $\Phi: (L,u) \mapsto \eps(L,u)$ induces a bijection
\[
\defo(M,R\ra S) \iso \Hom_{\cd_S A} (M,M\ten_S \n[1]).
\]
\end{proposition}

Clearly, the bijection of the proposition is equivariant with
respect to the action of $\Aut_{\cd_S A}(M)$. Therefore we have the

\begin{corollary} The map $\Phi: (L,u)\mapsto \eps(L,u)$ induces a
bijection
\[
\defo'(M,R \ra S) \iso \Hom_{\cd_S A}(M,M\ten_S \n[1])/\Aut_{\cd_S A}(M).
\]
\end{corollary}

\begin{proof} We construct a map $\Psi$ which will turn out to be the inverse bijection.
We may and will assume that $M$ is $S$-relatively cofibrant in the category of dg $A$-modules,
i.e. it satisfies the $S$-relative variant of property (P) of \cite[3.1]{Keller94}.
This means \cite[7.5]{Keller98}
that $M$ admits an increasing filtration by dg $A$-submodules $M_n$, $n\geq 0$,
such that each inclusion $M_n \subset M_{n+1}$ splits as a morphism of graded
$A$-modules and the subquotient $M_{n+1}/M_n$ is isomorphic to a direct summand
of a module $K\ten_S A$, where $K$ is a dg $S$-module.
Let a morphism
\[
f: M \ra M\ten_S\n[1]
\]
of $\cd_S A$
be given. Since $M$ is $S$-relatively cofibrant as a dg $A$-module, the morphism
$f$ is realized by a map $\tilde{f}$ of dg $A$-modules. We put $L=M\ten_S R= (M\ten_S\n) \oplus M$
as a graded $A\ten_S R$-module and we define its differential by
\[
d_L = \left[ \begin{array}{cc} d_{M\ten_S \n} & f \\ 0 & d_\n \end{array} \right].
\]
Clearly $L$ is free as a graded $R$-module (since $M$ is free as a graded $S$-module)
and we have an obvious isomorphism $u: L\ten_R S \ra M$.
So we have constructed an object $(L,\phi)$ of $\cf$
by choosing a representative $\tilde{f}$ of the homotopy class $f: M \ra M\ten_S\n[1]$. Let
us check that the connected component of $(L,\phi)$ is independent of the choice of the
representative. Let $\tilde{f}'$ be another choice and let $h: M \ra M\ten_S\n$ be a morphism
of graded $A$-modules such that $f' = f + d \circ h + h \circ d$. Let
$(L',u')$ be the object of $\cf$ constructed from $f'$. Then $L'$ equals $L$ as a graded
$A\ten_S R$-module and both equal $M\ten_S \n \oplus M$. The matrix
\[
\left[\begin{array}{cc} \id_{M\ten_S\n} & h \\ 0   & \id_M \end{array} \right]
\]
defines an isomorphism $v$ of dg $A\ten_S R$-modules from $L$ to $L'$ and we clearly have
$u' = u\circ (v\ten_R S)$. By definition, the map $\Psi$ sends $f$ to the connected component
of $(L,u)$. The easy check that $\Phi \circ \Psi$ is the identity is left to the
reader. To finish the proof, it is enough to check that $\Psi$ is surjective. For this,
let an object $(L',u')$ of $\cf$ be given. We may and will assume that $L'$ is
$R$-relatively cofibrant. Put  $M'=L'\ten_R S$.
Since $L'$ is free as a graded $R$-module, there is an isomorphism of
graded $R$-modules $M'\ten_S R \iso L'$ which lifts the identity of $M'$.
The differential of $L'$ then yields a differential of $M'\ten_S R = (M'\ten_S \n)\oplus M'$
given by a matrix
\[
\left[\begin{array}{cc} d_{M'\ten_S\n} & f' \\ 0 & d_{M'} \end{array} \right].
\]
Now define $f: M \ra M\ten_S\n[1]$ by the commutative square of $\cd_S A$
\[
\xymatrix{
M' \ar[d]_{u'} \ar[r]^(0.4){f'} & M'\ten_S \n[1] \ar[d]^{u'\ten_S \n[1]} \\
M              \ar[r]^(0.4){f}  & M \ten_S \n[1]
}
\]
Since $M'$ is $S$-relatively cofibrant, there is a morphism of dg modules $u_1: M' \ra M$
lifting $u'$. Moreover, $f\circ u_1$ is homotopy equivalent to
$(u_1\ten_S\n[1])\circ f'$. Choose an homotopy $h$ between the two. Then the matrix
\[
\left[\begin{array}{cc} u_1\ten \id_\n & h \\ 0 & u_1 \end{array} \right]
\]
defines a map
\[
\tilde{u}: L'=(M'\ten_S \n) \oplus M' \ra (M\ten_S \n) \oplus M = L.
\]
This is in fact a morphism of dg $A\ten_S R$-modules and clearly it gives
a morphism of $\cf$. Since the triangle (\ref{CharTriangle}) does not exist
in the category $\cd_R A$, it is not immediate that $\tilde{u}$ is invertible.
However, starting from an inverse of $u'$, we can analogously construct
a morphism $\tilde{v} : L \ra L'$ of $\cd_R A$. Then the reduction of
$\tilde{u} \tilde{v}$ modulo $\n$ is homotopic to the identity. Let
$H$ be an homotopy. The homotopy $H$ is a morphism of graded $A$-modules, and we
can lift it to a morphism of graded $A\ten_S R$-modules. We see that $\tilde{u}\tilde{v}$
is homotopic to a morphism $w$ whose reduction modulo $\n$ is the identity.
Since $L'$ and $L$ are free over $R$ and $\n$ is nilpotent, it follows that $w$ is invertible.
So $\tilde{u}\tilde{v}$ is homotopic to an invertible morphism. Thus $\tilde{u}\tilde{v}$ is
invertible in the homotopy category. Similarly, we see that $\tilde{v}\tilde{u}$
is invertible in the homotopy category.

We conclude that $(L',u')$ is in the isomorphism class
of $\cf$ which is the image of $f$ under $\Psi$. Hence $\Psi$ is surjective.
\end{proof}

\mysection{An exact sequence}{AnExactSequence}
Let $T$ be a commutative dg algebra, $S$ an augmented dg $T$-algebra and $R=S\oplus \n$
an augmented dg $S$-algebra. Thus $R$ becomes an augmented dg $T$-algebra $R=T\oplus \m$
and we have a commutative diagram with exact rows
\[
\xymatrix{
0 \ar[r] & \n \ar[d] \ar[r]     & R \ar@{=}[d] \ar[r] & S \ar[r] \ar[d]  & 0 \\
0 \ar[r] & \m \ar[r]            & R            \ar[r] & T \ar[r] & 0 .
}
\]
Let $A$ be a dg $T$-algebra free as a graded $T$-module and let $M$ be
a dg $A$-module free as a graded $T$-module.
Then the tensor product $?\ten_R S$ yields a natural map
\[
\defo(M, R\ra T) \ra \defo(M, S\ra T).
\]
On the other hand, if we have a representative $(L',u')$ of an element
of $\defo(M\ten_T S, R \ra S)$, we obtain an element of $\defo(M,R\ra T)$
by taking $L=L'$ and $u$ the composition
\[
L\ten_S T = L'\ten_S T \longarr{u'\ten_S T} (M\ten_T S)\ten_S T = M.
\]
The following lemma is immediate from these definitions.
\begin{lemma}
\begin{itemize}
\item[a)] The sequence
\[
\defo(M\ten_T S, R \ra S) \ra \defo(M, R\ra T) \ra \defo(M, S\ra T) \ra *
\]
is exact in the sense that the second map is surjective and its fibre over the
base point is the image of the first map.
\item[b)] The sequence
\[
* \ra \defo'(M\ten_T S, R \ra S) \ra \defo'(M, R\ra T) \ra  \defo'(M,S\ra T) \ra *
\]
is exact in the sense that the second map is surjective, and its fiber over the
base point equals the image of the first map, which is injective.
\end{itemize}
\end{lemma}

\mysection{A base change isomorphism}{BaseChangeIsomorphism}
Let $T$ be a commutative dg algebra, $S$ an augmented $T$-algebra and $R=S\oplus \n$
an augmented $S$-algebra. Thus $R$ becomes an augmented $T$-algebra $R=T\oplus \m$.
Let $R'=T\oplus \n$ so that we have a commutative diagram with exact rows
\[
\xymatrix{
0 \ar[r] & \n \ar@{=}[d] \ar[r] & R            \ar[r] & S \ar[r] & 0 \\
0 \ar[r] & \n \ar[r]            & R' \ar[u]    \ar[r] & T \ar[u] \ar[r] & 0 .
}
\]
Let $A$ be a dg $T$-algebra free as a graded $T$-module and let $M$ be
a dg $A$-module free as a graded $T$-module. The tensor products $?\ten_{R'} R$
and $?\ten_T S$ yield a natural map
\[
\defo(M, R' \ra T) \ra \defo(M\ten_T S, R \ra S).
\]
\begin{lemma}
If $\n^2=0$, this map is an isomorphism.
\end{lemma}

\begin{proof} We have a commutative square
\[
\begin{array}{ccc}
\defo(M, R'\ra T) & \ra & \defo(M\ten_T S, R\ra S) \\
\da       &     & \da \\
\Hom_{\cd_T A}(M, M\ten_T \n[1]) & \ra & \Hom_{\cd_S A}(M\ten_T S, M\ten_T S \ten_S \n[1])
\end{array}
\]
whose vertical maps are given by proposition (\ref{InfDefModules}) and whose horizontal
maps are given by the tensor functor. So the vertical maps are bijective.
And the lower horizontal arrow identifies
with the adjunction isomorphism
\[
\Hom_{\cd_T A}(M, M\ten_T \n[1])  \ra \Hom_{\cd_S A}(M\ten_T S, M\ten_T \n[1])
\]
where we view $M\ten_T \n[1]$ on the left hand side as the restriction to
$A$ of the $A\ten_T S$-module $M\ten_T \n[1]$ on the right hand side.
\end{proof}

\begin{corollary} With the above notations, suppose that the map
\[
\defo(M\ten_T S, R \ra S) \ra \defo'(M\ten_T S, R \ra S)
\]
is bijective. Then we have an exact sequence
\[
* \ra \defo(M, R' \ra T) \ra \defo'(M, R \ra T) \ra \defo'(M, S \ra T) \ra *
\]
in the sense that the second map is surjective and its fibre over the base point
is the image of the first map, which is injective.
\end{corollary}

This follows from the lemma and from part b) of
lemma~\ref{AnExactSequence}.

\mysection{Infinitesimal deformations of the bimodule $A$}{InfDefA}

Let $A$ be a dg $k$-algebra, $T$ a commutative dg $k$-algebra and $S=T\oplus \n$
an augmented commutative dg $T$-algebra with $\n^2=0$.

We consider the module $M=A\ten_k T$ over the algebra $B=A\op\ten A \ten T$.

\begin{lemma}Suppose that $\Hs T$ is of finite total dimension. Then
the canonical map
\[
\defo(M, S \ra T) \ra \defo'(M, S \ra T)
\]
is bijective.
\end{lemma}

\begin{proof} By (\ref{InfDefModules}), we have to show that the action
of $\Aut_{\cd_T B}(A\ten T)$ on
\[
\Hom_{\cd B}(M, M\ten_T \n[1])
\]
is trivial.
We claim that we have an isomorphism
\[
\R\Hom_{B} (A \ten T, A \ten T) \ten_T \n[1] \iso
\R\Hom_{B} (A\ten T, A \ten T \ten_T \n[1]).
\]
Indeed, we have an adjunction isomorphism
\[
\R\Hom_{B} (A \ten T, A \ten T) \ten_T \n[1] \iso
\R\Hom_{A^e}(A, A\ten T) \ten_T \n[1].
\]
Now since $\Hs T$ is of finite total dimension, we have an isomorphism
\[
\R\Hom_{A^e}(A, A\ten T)  \iso \R\Hom_{A^e}(A,A)\ten T.
\]
Combining the two preceding isomorphisms, we obtain an isomorphism
\[
\R\Hom_{B} (A \ten T, A \ten T) \ten_T \n[1] \iso
\R\Hom_{A^e}(A , A) \ten \n[1].
\]
Now we have isomorphisms
\[
\R\Hom_{A^e}(A , A) \ten \n[1] \iso \R\Hom_{A^e}(A,A\ten\n[1]) \iso
\R\Hom_{B}(A\ten T, A\ten T\ten_T \n[1]) \ko
\]
where we have first used that $\Hs \n$ is of finite total dimension and
then the adjunction, as above. The claim follows.
Now $\Aut_{\cd_T B}(A\ten T)$ is the group of invertible elements
of the zeroth homology of
\[
\R\Hom_{B} (A \ten T, A \ten T).
\]
The claim follows since this dg algebra is commutative up to homotopy.
\end{proof}

\mysection{Definition of the Lie bracket}{DefLieBracket}
Let $A$ be a dg $k$-algebra. Let $\m$ be a dg $k$-module whose homology
is of finite total dimension. Let $R=k\oplus\m$ denote the augmented commutative
dg algebra with $\m^2=0$. We consider $A$ as an $A$-$A$-bimodule. We define
\[
G(\m) = \defo'(A, k\oplus \m \ra k).
\]
In other words, $G(\m)$ is formed by the isomorphism classes of
objects in $\cd_R(A\op\ten A \ten R)$ whose reduction modulo $\m$
is isomorphic to $A$ in $\cd(A\op\ten A)$.
Note that, according to lemma~\ref{InfDefA}, we have a canonical bijection
\[
\defo(A, k\oplus\m \ra k) \iso \defo'(A, k\oplus\m \ra k).
\]
We will view $\cd(A\op\ten A)$ (resp.~$\cd_R(A\op\ten A\ten R)$) as a monoidal
category for the derived tensor product over $A$ (resp.~for the relative derived
tensor product over $A\ten R$). The monoidal structure of
$\cd_R(A\op\ten A \ten R)$ induces a monoid structure on $G(\m)$ and the bijection
\[
G(\m) \liso \defo(A, k\oplus\m \ra k) \ra \Hom_{\cd(A\op\ten A)}(A, A\ten \m[1])
\]
is a monoid morphism. In particular, $G(\m)$ is an abelian group, functorial in $\m$.

For two dg $k$-modules $\m_1$ and $\m_2$ whose homology is of finite
total dimension, we define a Lie bracket
\[
G(\m_1) \times G(\m_2) \ra G(\m_1 \ten \m_2)
\]
as follows : Let $L_1$ and $L_2$ represent elements of $G(\m_1)$
resp. $G(\m_2)$. Put $R_i=k\oplus \m_i$. Let $U_i$ be the image of
$L_i$ in $\cd_R (R\ten A\op\ten A)$ where $R=R_1\ten R_2$
(note that the kernel of $R\ra k$ is not of square zero !).
The $U_i$ are invertible objects of a monoidal category. Let
$V$ denote the commutator of $U_1$ with $U_2$. Then $V$ yields
an element of $\defo'(A, R \ra k)$. We have a
canonical map
\[
G(\m_1 \ten \m_2) \ra \defo'(A, R \ra k)
\]
and we claim that it is injective and that $V$ lies in its image.
Indeed, the image of $V$ in $\defo'(A, R_1 \ra k)=G(\m_1)$ vanishes
since $G(\m_1)$ is commutative. Thus $V$ lies in
\[
G(\m_2\oplus \m_1\ten \m_2)=\defo'(A, k\oplus (\m_2 \oplus \m_1 \ten \m_2))
\]
by (\ref{AnExactSequence}). The image of $V$ in $G(\m_2)=\defo'(A, k\oplus \m_2\ra k)$
also vanishes since $G(\m_2)$ is commutative. So again by~\ref{AnExactSequence},
$V$ lies in $\defo'(A, k\oplus (\m_1\ten\m_2) \ra k)=G(\m_1\ten \m_2)$.

\mysection{From coalgebra automorphisms to bimodules}{CoalgAutomBimod}
Let $R$ be a commutative dg $k$-algebra and $A$ a (typically noncommutative) dg $R$-algebra.
Denote by $SA$ the graded $R$-module with $(SA)^p=A^{p+1}$. We recall the {\em bar
construction} of $A$ relative to $R$. It is the dg $R$-coalgebra $C^+$ defined as follows:
Its underlying graded space is
\[
R \oplus SA \oplus (SA\ten_R SA) \oplus \ldots \oplus (SA)^{\otimes_R p} \oplus \ldots \;\; .
\]
The comultiplication is defined by
\begin{eqnarray*}
\Delta(a_1, \ldots, a_p) & = & 1_R \ten (a_1, \ldots, a_p) +
                 \sum_{i=1}^{p-1} (a_1, \ldots, a_i)\ten (a_{i+1}, \ldots, a_p)  \\
                          &  &  + (a_1, \ldots, a_p)\ten 1_R \;\; .
\end{eqnarray*}
Moreover $C^+$ is endowed with the counit $\eta: C^+ \ra R$ given by the natural projection and
the coaugmentation $\eps: R \ra C^+$ given by the natural inclusion. The composition of the
projection $C^+ \ra SA$ with the canonical morphism $s: SA \ra A$ of degree $+1$ yields
a homogeneous morphism $\tau : C^+ \ra A$ of degree $+1$. A {\em coderivation}
of $C^+$ is a homogeneous
$R$-linear map $D: C^+ \ra C^+$ such that $\Delta \circ D = D\ten \id + \id \ten D$.
Note that this implies that $\eta\circ D =0$. Let $\Coder(C^+,C^+)$ denote the graded
$R$-module whose $p$-th component is formed by the
coderivations of degree $p$. Then the composition with $\tau$ is a bijection
onto the space of homogeneous $R$-linear morphisms from $C^+$ to $A$
\begin{equation}
\Coder(C^+,C^+) \iso \Hom_R(C^+, A) \ko D \mapsto \tau\circ D. \label{CoderGerst}
\end{equation}
In particular, $C^+$ admits a unique coderivation $d_{C^+}$ of degree $+1$ such that
$\tau \circ d_{C^+}$ restricted to $(SA)^{\ten_R p}$ vanishes for $p\neq 1,2$, and equals
\[
\mu \circ (s \ten s) : (SA)\ten_R (SA) \ra A
\]
for $p=2$ and $-d_A\circ s$ for $p=1$. Here $\mu$ denotes the multiplication
of $A$. The facts that $\mu: A\ten A \ra A$ is a morphism of complexes and
that $\mu$ is associative are equivalent to the fact that $d_{C^+}^2=0$.

\begin{proposition}[\cite{Stasheff93}]
Endowed with the supercommutator and the differential $D \mapsto [d_{C^+},D]$
the graded space $\Coder(C^+,C^+)$ becomes a differential graded Lie algebra
which is isomorphic to the Gerstenhaber Lie algebra by the map
(\ref{CoderGerst}).
\end{proposition}

For two homogeneous $R$-linear morphisms $f,g: C^+ \ra A$, we define $
f\star g = \mu \circ (f\ten g) \circ \Delta$.
Then we have
\begin{equation}
\tau \star \tau = d\circ \tau + \tau \circ d \quad\quad  \label{twisting}
\end{equation}
For a dg right $A$-module $M$ and a dg left $C^+$-comodule $N$, we denote by $M\ten_\tau N$
the dg $R$-module $M\ten_R N$ with the differential defined by
\[
d(x\ten y) = d(x) \ten y + (-1)^p x \ten y + \sum x \,\tau(y_{(1)}) \ten y_{(2)} \ko
\]
where $x$ is homogeneous of degree $p$ and $\delta(y)=\sum y_{(1)}\ten y_{(2)}$
(Sweedler's notation). Similarly, for a left $A$-module $M$ and a right $C^+$-comodule $N$,
we define $N\ten_\tau M$ to be $N\ten_R M$ with the differential defined by
\[
d(x\ten y) = d(x) \ten y + (-1)^p x \ten y - \sum x_{(1)} \ten \tau(x_{(2)})\, y  \ko
\]
(note the sign in front of $\sum$). The fact that the squares of these differentials vanish follows
from equation (\ref{twisting}).

The dg $R$-module $A\ten_\tau C^+$ inherits a right $C^+$-comodule structure from
$C^+$ and a left $A$-module structure from $A$. It yields the dg $A$-$A$-bimodule
\[
A \ten_\tau C^+ \ten_\tau A.
\]
It is not hard to check that up to the signs of the differentials, this is
the (sum) total dg module associated with the bar resolution of the $A$-$A$-bimodule $A$.
In particular, we have a canonical quasi-isomorphism (which is even an homotopy equivalence
of left dg $A$-modules or right dg $A$-modules)
\[
A \la A \ten_\tau C^+ \ten_\tau A.
\]

Now let $\phi: C^+ \ra C^+$ be an automorphism of the dg counital
$R$-coalgebra $C^+$. Define $C^+_\phi$ to be the dg $C^+$-$C^+$-bicomodule whose
left comultiplication is that of $C^+$ whereas the right comultiplication is
$(\id\ten \phi) \circ \Delta$. We define the bimodule
\[
X(\phi) = A \ten_\tau  C^+_\phi \ten_\tau A.
\]
Note that the underlying graded module of $X(\phi)$ is $A\ten_R C^+ \ten_R A$ but
the differential is twisted by $\phi$. Now let $\psi: C^+ \ra C^+$ be another
automorphism. Then we have a natural morphism of dg $C^+$-$C^+$-bicomodules
\[
 C^+_{\psi \phi} \ra  C^+_\phi\ten_\tau A \ten_\tau  C^+_\psi \ko
c \mapsto   \sum c_{(1)} \ten 1_A \ten \phi(c_{(2)}).
\]
It induces a morphism of dg $A$-$A$-bimodules
\begin{equation}
X(\psi \phi) \ra X(\phi) \ten_A X(\psi). \label{GroupHom}
\end{equation}

\begin{proposition} \begin{itemize}
\item[a)] As a left dg $A$-module, $X(\psi)$ is relatively cofibrant. In particular, we have a canonical isomorphism
\[
X(\phi) \lrelten_A X(\psi) \iso X(\phi) \ten_A X(\psi)
\]
in $\cd_R(A\op\ten A)$.
\item[b)] The morphism (\ref{GroupHom}) is a homotopy equivalence of dg $R$-modules.
\end{itemize}
\end{proposition}
In the next section, we will deduce this from results of
\cite{Lefevre2002}. Note that we obtain a morphism of groups
from the group
\[
\Aut_{R-\mbox{\scriptsize coalg}}(C^+)\op
\]
to the group of autoequivalences of the relative derived category
$\cd_R A$.

Now suppose that $R$ is an augmented dg $k$-algebra and
$R=k\oplus \m$ the corresponding decomposition. Suppose
that $\phi\ten_R k : C^+\ten_R k \ra C^+ \ten_R k$ is
the identity. Then clearly $X(\phi)\ten_R k$ is isomorphic
to $X(\id)\ten_R k$ as a dg $A$-$A$-bimodule and we have a
canonical isomorphism $X(\id)\ten_R k \ra A$ in $\cd(A\op\ten A)$.
So we obtain a canonical isomorphism $u_\phi :X(\phi) \ra A$
in $\cd(A\op\ten A)$ and an object $(X(\phi), u_{\phi})$
of the fiber category $\cf$ associated with the reduction
functor $\cd_R(R\ten A\op\ten A) \ra \cd(A\op\ten A)$
(\conf section \ref{InfDefModules}).

\mysection{Modules and comodules}{ModulesAndComodules}
Let $R$ be a commutative dg $k$-algebra
and $A$ a dg $R$-algebra of the form $A=A'\ten_k R$
for some dg $k$-algebra $A'$.
We define $A^+$ to be the augmented algebra $R\oplus A$ and
$C^+$ to be the coaugmented coalgebra defined in the previous
section. We still denote by $\tau : C^+ \ra A^+$ the composition
of the morphism $\tau$ of the previous section with the
inclusion $A \ra A^+$. Denote by $\Comod C^+$ the category of dg
counital right $C^+$-comodules.
% These comodules are automatically cocomplete since $C^+$ is cocomplete !
\[
\xymatrix{
\Mod A^+ \ar@<1ex>[d]^B \\
\Comod C^+ \ar@<1ex>[u]^\Omega
}
\]
given by
\[
BM = M\ten_\tau C^+ \ko \Omega N = N\ten_\tau A^+ .
\]
One can check \cite{Lefevre2002} that they form an adjoint pair.
Let $\cw$ denote the class of morphisms $s$ in $\Mod A^+$ whose restriction
to $\Mod R$ is an homotopy equivalence and let
$\cw'$ be the class of morphisms $s$ of $\Comod C^+$ such that $\Omega s$
belongs to $\cw$.

\begin{theorem}
\begin{itemize}
\item[a)] The dg $A$-module $\Omega N$ is relatively cofibrant for each
dg $C^+$-comodule $N$.
\item[b)]
We have $B \cw\subset \cw'$ and $\Omega \cw' \subset \cw$ and the functors
$B$ and $\Omega$ induce quasi-inverse equivalences between the localized categories
\[
(\Mod A^+)[\cw^{-1}] \iso (\Comod C^+)[{\cw'}^{-1}].
\]
\end{itemize}
\end{theorem}

\begin{theorem}
The restriction functor $\Mod A \ra \Mod A^+$ induces an equivalence
from $\cd_R(A)$ onto the full subcategory of $(\Mod A^+)[\cw^{-1}]$
whose objects are the dg modules $M$ such that $M\ten_\tau C^+$ is
$R$-relatively acyclic (i.e. its underlying dg $R$-module is contractible).
\end{theorem}

These theorems are proved in \cite{Lefevre2002} in the case where $R=k$
(the first one corresponds to Theorem~2.2.2.2 and the second one to
Proposition~4.1.2.10 in \cite{Lefevre2002}).
We omit the proof in the general case since it is similar. Note however
the following: If $M$ is a dg $A$-module, then according to the first
theorem, we have a canonical $R$-relative quasi-isomorphism
\[
M\ten_\tau C^+ \ten_\tau A^+ \ra M.
\]
The existence of such a quasi-isomorphism is not surprising. Indeed,
the decomposition $A^+=A\oplus R$ yields an $R$-split short exact sequence
\[
0\ra M\ten_\tau C^+ \ten_\tau A \ra M\ten_\tau C^+ \ten_\tau A^+
\ra M\ten_\tau C^+ \ten R \ra 0.
\]
The last term identifies with the augmented bar resolution of $M$. It is therefore
relatively $R$-acyclic. The first term identifies with the bar resolution of $M$.
It is therefore relatively $R$-quasi-isomorphic to $M$.

Now suppose that $\phi: C^+ \ra C^+$ is an automorphism of dg $R$-coalgebras. Then it induces
a selfequivalence $F_\phi$ of $\Comod C^+$ given by twisting by $\phi$, i.e. if $N$ is a
dg $C^+$-comodule, then $F_\phi N$ is the dg comodule with the same underlying graded module and the
same differential but with the new comultiplication $\delta_\phi = (\id \ten \phi)\circ \delta$.
Clearly $F_\phi$ preserves the subcategory of $R$-relatively acyclic comodules. Thus
the composition $\Omega F_\phi B$ preserves the image of $\cd_R(A)$ in $(\Mod A^+)[\cw^{-1}]$.
More precisely, if $M$ is a dg $A^+$-module then
\[
\Omega F_\phi B M = M\ten_\tau C^+_\phi \ten_\tau A^+
\]
and if $M$ comes from a (unital) dg $A$-module, then the last term is $R$-relatively quasi-isomorphic
to its submodule $M\ten_\tau C^+_\phi \ten_\tau A$ since
\[
M\ten_\tau C^+_\phi \ten_\tau R = M\ten_\tau C^+
\]
is $R$-relatively acyclic. So for each $M\in \Mod A$, we have a canonical $R$-relative
quasi-isomorphism
\[
M\ten_\tau C^+_\phi \ten_\tau A \ra \Omega F_\phi B M.
\]
If $\psi$ is another automorphism, by composition,
we obtain a canonical $R$-relative quasi-isomorphism
\[
M\ten_\tau C^+_\phi \ten_\tau A \ten_\tau C^+_\psi \ten_\tau A
\ra \Omega F_\psi B \Omega F_\phi B M.
\]
On the other hand, the adjunction morphism $\id \ra B\Omega$
yields a morphism
\[
\Omega F_\psi B \Omega F_\phi B M \la \Omega F_\psi F_\phi B M
= \Omega F_{\psi \phi} B M \ko
\]
which is also an $R$-relative quasi-isomorphism. To prove
proposition~\ref{CoalgAutomBimod}, it remains to
be checked that the morphism (\ref{GroupHom}) constructed in (\ref{CoalgAutomBimod})
makes the following square commutative
\[
\xymatrix{
A\ten_\tau C^+_\phi \ten_\tau A \ten_\tau C^+_\psi \ten_\tau A  \ar[r] &
        \Omega F_\psi B \Omega F_\phi B A \\
A\ten_\tau C^+_{\psi\phi} \ten_\tau A \ar[u] \ar[r]  &
        \Omega F_{\psi \phi} B A \ar[u]
}
\]
This is left to the reader.

\mysection{Proof of the main result}{ProofMainResult}
Let $A$ be a dg $k$-algebra and $R$ a commutative augmented dg $k$-algebra.
We write $R=k\oplus \m$, where $\m$ is the kernel of the augmentation.

Let $C^+$ be the bar construction of $A$ relative to $k$ (\conf \ref{CoalgAutomBimod}).
Then the bar construction of $A\ten R$ relative to $R$ identifies with $C^+\ten R$.
We put
\[
\lieaut(C^+,R) = \ker (\Aut_R(C^+\ten R) \ra \Aut_k(C^+))
\]
where $\Aut_R$ denotes the group of automorphisms of dg counital $R$-coalgebras.

We define $\liedpic(A,R)$ to be the group of isomorphism classes
(\conf section~\ref{MainResult}) of invertible dg bimodules $X\in
\cd_R(A\op\ten A\ten R)$ free as graded $R$-modules such that $X\ten_R
k$ is isomorphic to $A$ in $\cd(A\op\ten A)$.

By section~\ref{CoalgAutomBimod}, we obtain a morphism of groups
\[
\Phi  : \lieaut(C^+,R)\op \ra \liedpic(A,R)
\]
which is clearly functorial in $R$. As in the case of algebraic
groups, one obtains canonical Lie
brackets on the restrictions of these functors
to the subcategory of augmented dg $k$-algebras $R=k\oplus\m$ with
$\m^2=0$ (\conf section~\ref{DefLieBracket})
and $\Phi$ is compatible with the bracket.

\begin{lemma}If $R=k\oplus \m$ with $\m^2=0$, there is natural isomorphism
\[
\Zy{0} (\Coder_k(C^+, C^+)\ten_k \m) \ra \lieaut(C^+,R).
\]
\end{lemma}

This is a variant of a classical result on infinitesimal
deformations. The easy proof is left to the reader. As we recalled
from \cite{Stasheff93} in section \ref{CoalgAutomBimod}, we have a
natural isomorphism of dg Lie algebras
\[
C^{\cdot+1}(A,A) \iso \Coder_k(C^+,C^+).
\]
So for $\m^2=0$, we obtain morphisms
\[
\Zy{0} (C^{\cdot+1}(A,A)\ten \m) \ra \lieaut(C^+,R) \ra \liedpic(A,R)\op.
\]
compatible with the bracket. Now by \ref{InfDefModules}, we have an isomorphism
\[
\Phi : \liedpic(A,R) \iso \Hom_{\cd(A\op\ten A)}(A, A\ten\m[1]).
\]
It is easy to see that the composition
\[
\Zy{0} (C^{\cdot+1}(A,A)\ten \m) \ra  \liedpic(A,R) \iso  \H{0} (C^{\cdot+1}(A,A)\ten \m) .
\]
is the canonical surjection. So we have a commutative square
\[
\xymatrix{
\Zy{0}(C^{\cdot+1}(A,A)\ten \m) \ar[d] \ar[r] & \liedpic(A,R) \ar[d]^{\sim} \\
\H{0}(C^{\cdot+1}(A,A)\ten\m) \ar[r]^(0.45){\sim}   & \Hom_{\cd(A\op\ten A)}(A, A\ten\m[1])
}
\]
We see that if we transport the Gerstenhaber bracket to the lower right hand corner,
then the map
\[
\Phi : \liedpic(A,R)\op \ra \Hom_{\cd(A\op\ten A)}(A,A\ten\m[1])
\]
becomes an isomorphism which respects the bracket and is functorial with respect
to $R$ and with respect to invertible bimodule complexes $X\in\cd(A\op\ten B)$.

\end{document}